\documentclass[12pt,letterpaper]{article}
\usepackage[hmargin=1in,vmargin=1in]{geometry}

\usepackage{amsmath}
\usepackage{amsfonts}
\usepackage{amsthm}
\usepackage{amssymb}
\usepackage{array}
\usepackage{caption}
\usepackage{color}
\usepackage{graphicx}
\usepackage{hyperref}
\usepackage{verbatim}

\newtheorem{prop}{Proposition}

\newtheorem{hint}[prop]{Hint}

\usepackage{rotating}

\begin{document}

\title{Thinking Inside and Outside the Box}
\author{Tanya Khovanova}
\maketitle

\begin{abstract}
I discuss puzzles that require thinking outside the box. I also discuss the box inside of which many people think.
\end{abstract}

\section{Nine Dots Puzzle}
The following puzzle started the expression \textit{inside the box} \cite{K}.

\begin{quote}
Connect the dots by drawing four straight, continuous lines that pass through each of the nine dots without lifting the pencil from the paper.
\end{quote}

\begin{center}
\includegraphics[scale=0.2]{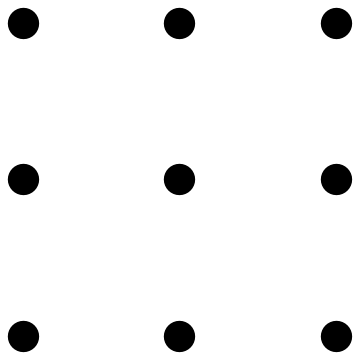}
\end{center}

Most people try something like this:

\begin{center}
\includegraphics[scale=0.3]{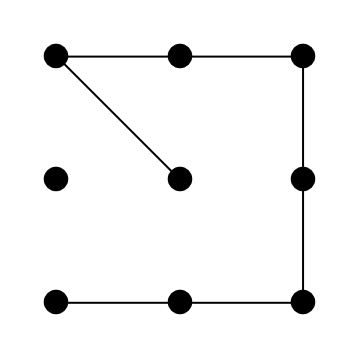}
\end{center}

and they fail to connect all the dots. They try to connect the dots that are formed by line segments that fit inside the square box around the dots:

\begin{center}
\includegraphics[scale=0.3]{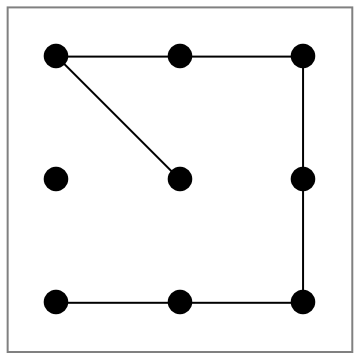}
\end{center}

They mentally restrict themselves to solutions that are literally inside this box. In the correct solution the line segments should be drawn outside this imaginary box:

\begin{center}
\includegraphics[scale=0.3]{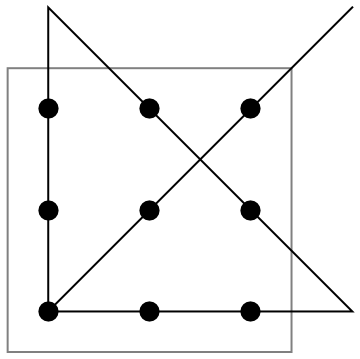}
\end{center}

There are many puzzles where people make assumptions that are not required by the puzzle. They constrain themselves to the \textit{inside of the box}. Let us look at such puzzles.

\section{Other Example of the Outside-the-Box Puzzles}

These puzzles have a common feature. Many people who try to solve them assume something that is not stated in the puzzle and therefore they fail to find the solution. 

\subsection{Other Nine Dots}

There is another puzzle with the same nine-dots setup.

\begin{quote}
What is the smallest number of squares needed to make sure that each dot is in its own region?
\end{quote}

\begin{center}
\includegraphics[scale=0.2]{9dots2.png}
\end{center}

Usually people who try to solve this puzzle come up with the following solution with four squares.

\begin{center}
\includegraphics[scale=0.3]{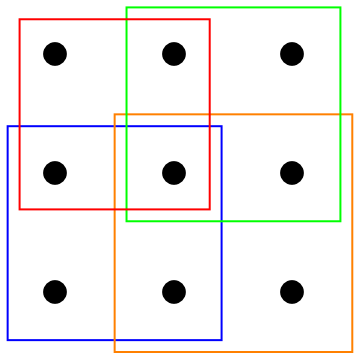}
\end{center}

As with the previous puzzle they imagine the dots are on a grid and try to build squares that have sides parallel to the grid lines. What is the outside-the-box idea? The sides of squares do not need to be parallel to the grid. This way we can find a solution with three squares.

\begin{center}
\includegraphics[scale=0.3]{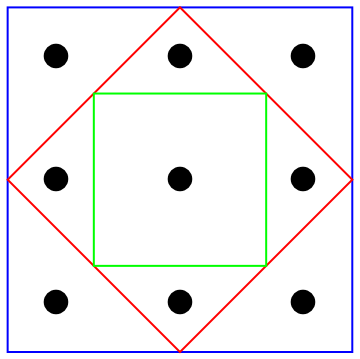}
\end{center}

\subsection{An Interesting Bet}

Here is another puzzle:

\begin{quote}
``I will bet you \$1'' said Fred, ``that if you give me \$2, I will
give you \$3 in return.'' Tom agreed and gave Fred \$2. How much money did Tom win?
\end{quote}

It seems that Fred will give Tom \$3 dollars at which point Tom will lose the bet and will have to give \$1 back. Overall Tom wins nothing. But there is a box around this puzzle. Most people assume that as Fred made the bet, he will follow on it. Actually, if Fred loses the bet he wins \$1 and Tom loses \$1.

\subsection{The Ancient Outside-the-Box Puzzle}

The following problem can be found in eighth-century writings.

\begin{quote}
A man has to take a wolf, a goat, and some cabbage across a river. His rowboat has enough room for the man plus either the wolf or the goat or the cabbage. If he takes the cabbage with him, the wolf will eat the goat. If he takes the wolf, the goat will eat the cabbage. Only when the man is present are the goat and the cabbage safe from their enemies. All the same, the man carries wolf, goat, and cabbage across the river. How?
\end{quote}

\begin{center}
\includegraphics[scale=0.5]{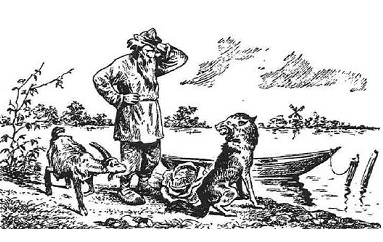}
\end{center}

I think that the reason the puzzle has survived for so many years is that the solution is based on an outside-the-box idea.

I don't want to provide a solution to this popular puzzle, but here is the non-trivial idea: Though the man wants to move his luggage to the other side, in one of the trips he needs to bring the goat back.

\subsection{Move a Digit}

Here is another puzzle:

\begin{quote}
In the equation $30 - 33 = 3$ move one digit to make it correct.
\end{quote}

The puzzle seems impossible. But the outside-the-box idea is to move a digit up to the exponent: $30 - 3^3 = 3$.

\subsection{Cigarette Butts}

Here is the last puzzle in this section:

\begin{quote}
A certain hobo who is skilled at making cigarettes can turn any 4 cigarette butts into a single cigarette. Today, this hobo has found 24 cigarette butts on the street. Assuming he smokes every cigarette he can, how many cigarettes will he smoke today? 
\end{quote}

On the surface, he can smoke $24/4 = 6$ cigarettes. What is the outside-the-box idea? He can reuse his own butts. After smoking 6 cigarettes, he will have 6 butts left. He can make one more cigarette. The answer is 7.

Or is it? What I love about this puzzle is that it has two layers. After smoking 7 cigarettes the hobo will have 3 butts left. There is another outside-the-box idea here. He can borrow a butt from a friend, smoke a cigarette and return the butt. At the end he can smoke 8 cigarettes.

\section{My Students}

I love outside-the-box puzzles. And I am good at them. Whenever I see such a puzzle, I always know the intended answer. Unfortunately, the moment I get the answer, I stop thinking about the puzzle. This is where my students come in. I always give such puzzles to my students, and they never fail to surprise me.

\subsection{A River-Crossing Puzzle}

Let's look at the following puzzle:

\begin{quote}
Two boys wish to cross a river, but there is a single boat that can take only one boy at a time. The boat cannot return on its own; there are no ropes or similar tricks; yet both boys manage to cross the river. How?
\end{quote}

The outside-the-box idea is that they started on different sides of the river. Many of my students do not see this answer. Nevertheless, they are very inventive and produce a lot of interesting different answers:

\begin{itemize}
\item There was another person on the other side of the river who brought the boat back.
\item There was a bridge.
\item The boys can swim.
\item They just wanted to cross the river and come back, so they did it in turns.
\end{itemize}

I gave a talk about thinking inside and outside the box at the 2016 Gathering for Gardner conference. I mentioned this puzzle and the inventiveness of my students. After my talk a guy approached me with another answer which is now my favorite:

\begin{itemize}
\item They wait until the river freezes over and walk to the other side.
\end{itemize}

\subsection{Apples in a Basket}

Here is another puzzle, over which my students didn't fail to amaze me.

\begin{quote}
You have a basket containing five apples. You have five hungry friends. You give each of your friends one apple.  After the distribution each of your friends has one apple each, yet there is an apple remaining in the basket. How can it be?
\end{quote}

The standard outside-the-box solution is to give an apple to a friend together with a basket. Here are some pearls from my students:

\begin{itemize}
\item Kill one of your friends.
\item You can't count.
\item You are narcissistic and you are one of your own friends.
\item You have two baskets, one has 5 apples, one has 1 apple.
\item One friend already has an apple.
\item An extra apple falls from a tree into the basket.
\item \textbf{My favorite:} the basket is your friend!
\end{itemize}

When I see these answers I regret that I stopped thinking about the puzzle and didn't see how many more ideas it can generate. At the same time I am elated at the inventiveness of my students. I wouldn't have considered killing one of the friends as an option, but then I'm not in the 8th grade like my student.

\subsection{The Original Nine-Dots puzzle}

After the G4G conference mentioned earlier, Jason Rosenhouse told me a solution for the original nine-dots puzzle that requires only three lines:

\begin{center}
\includegraphics[scale=0.3]{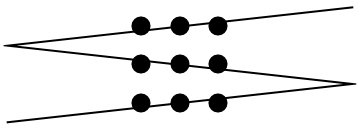}
\end{center}

The outside-the-box idea here is to use the thickness of the dots.

\subsection{An Irresistible Cannonball}

My students are inventive not only when solving outside-the-box puzzles. I am taking a detour here with another puzzle that is not the outside-the-box type. But I love what my student suggested and think you will too.

\begin{quote}
What happens if an irresistible cannonball hits an immovable post?
\end{quote}

This puzzle is known as the Irresistible Force Paradox. I borrowed it from the book \textit{What is the name of this book?} by Raymond Smullyan \cite{Smullyan}. The standard answer is that the given conditions are contradictory and the two objects cannot exist at the same time.

This is what one of my students wrote:

\begin{quote}
The post falls in love with the cannonball as it is so irresistible.
\end{quote}

\section{Where is the Box?}

I am good at thinking outside the box. I even drew a picture of myself to represent this.

\begin{center}
\includegraphics[scale=0.3]{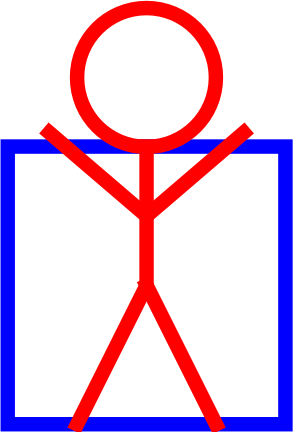}
\end{center}

Unfortunately I have to conclude, that I am inside my own, bigger box.

\begin{center}
\includegraphics[scale=0.3]{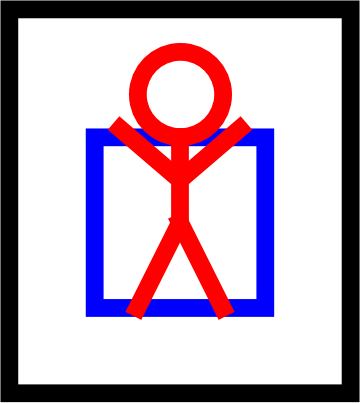}
\end{center}

All of us have our own boxes. It is good that we can learn from each other about the beauty outside our boxes.

\section{More Outside-the-Box puzzles}

There are many more outside-the-box puzzles. The fact that you know that they need an outside-the-box solution will help you find it. Actually each time you are stuck on a puzzle, it makes sense to assume that the puzzle might need an outside-the-box idea.

\begin{quote}
\textbf{Puzzle.} Four matchsticks form a square. How many non-overlapping squares can be formed using eight matchsticks?
The matchsticks do not intersect each other and they can't be broken.
\end{quote}

As you might have guessed the answer is more than 2. If you need a hint, you need a mirror.

\reflectbox{\rotatebox{180}{\textbf{Hint.} The squares are smaller than you might think.}}

\begin{quote}
\textbf{Puzzle.} I arranged eight sticks in the shape of a fish. What is the minimum number of sticks that must be moved to make the fish face the opposite direction?
\begin{center}
\includegraphics[scale=0.3]{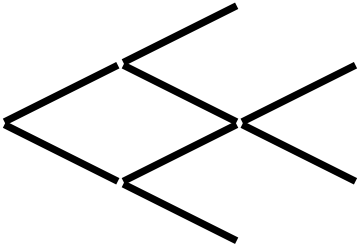}
\end{center}
\end{quote}

\reflectbox{\rotatebox{180}{\textbf{Hint.} The fish moves down.}}

\begin{quote}
\textbf{Puzzle.} These four sticks make a glass with a cherry in it. Can you move just 2 sticks so that the cherry is outside the glass?
\begin{center}
\includegraphics[scale=0.3]{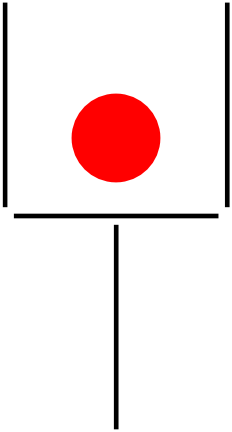}
\end{center}
\end{quote}

\reflectbox{\rotatebox{180}{\textbf{Hint.} The glass will be upside down.}}

\begin{quote}
\textbf{Puzzle.}
You have two wallets, each containing a quarter. Yet the total amount of money you have is 25 cents. How could this be?
\end{quote}

There are no hints for this one: it is too easy.


\begin{thebibliography}{9}

\bibitem{K} M.~Kihn, `Outside the Box': the Inside Story, FastCompany 1995.

\bibitem{Smullyan}
R.~Smullyan, {\it What is the Name of this Book?}, Prentice-Hall, 1978.



\end{thebibliography}
\end{document}